\documentclass[a4paper,12pt]{article}
\usepackage{amssymb,amsmath}
\newcommand{\Fb}{\overline{F}}
\newcommand{\Q}{\mathbb{Q}}

\newcommand{\A}{\mathbb{A}}

\newcommand{\cp}{{\mathfrak p}}
\newcommand{\Proj}{\mathbb{P}}
\renewcommand{\b}[1]{{\bf #1}}
\newcommand{\x}{\b{x}}
\newcommand{\bu}{\b{u}}

\newcommand{\Z}{\mathbb{Z}}
\newcommand{\lb}[1]{\mbox{\scriptsize\bf #1}}

\newcommand{\rank}{{\rm Rank}}

\newtheorem{theorem}{Theorem}
\newtheorem{corollary}{Corollary}

\newtheorem{lemma}{Lemma}

\begin{document}
\title{Zeros of Systems of $\cp$-adic Quadratic Forms}
\author{D.R. Heath-Brown\\Mathematical Institute, Oxford}
\date{}
\maketitle
\section{Introduction}

Let $K$ be a finite extension of $\Q_p$ with associated prime ideal
$\cp$, and let
$q^{(i)}[x_1,\ldots,x_n]\in K[x_1,\ldots,x_n]$ be quadratic forms, for
$1\le i\le r$.  It would follow from the conjecture of Artin
\cite[Preface]{Artin} that these forms have a simultaneous non-trivial zero in
$K^n$ providing only that $n>4r$.  Although Artin's conjecture is
known to be false in general (see Terjanian \cite{Ter}, for example),
this particular consequence of the conjecture is still open.  The cases
$r=1$ and $r=2$ have been successfully handled, the former being
due to Hasse \cite{Has} and the latter to Demyanov \cite{Dem}.  For $r=3$
it has been shown by Schuur \cite {Schuur} that $n\ge 13$ suffices
when the residue  field has odd characteristic and cardinality
at least 11.  No analogous result for $r\ge 4$ has been established
until now.  However it follows from the work of Ax and Kochen
\cite{AK} that if the degree $[K:\Q_p]=D$ is given, then $n\ge 4r+1$
variables suffice as soon as $p\ge p(r,D)$, for some prime $p(r,D)$.
The proof uses methods from mathematical logic, and does not yield a
practical value for $p(r,D)$.

If one is willing to allow more variables, then further results are
available.  Thus Martin \cite{mar} has shown that for any $K$ it is
sufficient that $n\ge 2r^2+3$ if $r$ is odd, and that $n\ge 2r^2+1$ if
$r$ is even.  One can do a little better for large $r$ but the bound
on $n$ is asymptotically $2r^2$ in all such results.

The purpose of the present paper is to develop an analytic method which
will establish the following result.
\begin{theorem}\label{main}
Let $K$ have residue  field $F$ and suppose that $\#F=q$.
Then the quadratic forms 
$q^{(1)},\ldots,q^{(r)}$ have a non-trivial common
zero over $K$ providing that $q\ge(2r)^r$.  More specifically it
suffices that $q> n\ge 4r+1$ and
$\sigma_1+\sigma_2<1$, where
\[\sigma_1=q^{r-n}+\sum_{\lceil n/2r\rceil-1\le t\le n/2}q^{-t}
(\frac{q}{2t+1})^{[4rt/n]}(2t+1)^r\]
and
\[\sigma_2=\frac{1}{q-1}\sum_{\rho=2(\lceil n/2r\rceil-1)}^{n-1}\,
\sum_{0\le t\le(n-\rho)/2}C_{\rho,t}q^{-\rho-t+[2r\rho/n]+[2r(\rho+2t)/n]}\]
with
\[C_{\rho,t}=(\rho+1)^{r-[2r\rho/n]}(2t+1)^{r-[2r(\rho+2t)/n]}.\]
\end{theorem}
Here we use the notation
\[\lceil \theta\rceil=\min\{n\in\Z:n\ge\theta\}.\]
Some small improvements in the values of $\sigma_1$ and $\sigma_2$ are
possible, but these have little effect on the range of $q$ which one
may handle.

It should be emphasized that the Ax-Kochen theorem gives no
information about fields with a fixed characteristic $p$.  Thus it
leaves open the possibility that Artin's conjecture is {\em never} true for
dyadic fields, for example. In contrast our result shows that it is
sufficient to have $\#F$ large enough.

We have the following corollary.  The case $r=8$ will be of relevance
later.
\begin{corollary}\label{c1}
It suffices to have
$n\ge 4r+1$ in the following cases.
\begin{enumerate}
\item[(i)] $r=3$ and $q\ge 37$; 
\item[(ii)] $r=4$ and $q\ge 191$;
\item[(iii)] $r=8$ and $q\ge 271919$;
\end{enumerate}.
\end{corollary}

As an indication of what can be achieved for larger values of $n$ we
investigate the condition $n>r^2$, which may be compared with Martin's result
\cite{mar} mentioned above in which one requires $n\ge 2r^2+3$ if $r$
is odd, and that $n\ge 2r^2+1$ if
$r$ is even, for any $q$.

\begin{corollary}\label{c2}
It suffices to have
$n\ge r^2+1$ providing that $r\ge 5$ and $q\ge (4\times 10^8)r^2$.
\end{corollary}

The coefficient in front of $r^2$ can certainly be improved, but the
importance of the result is that we require a lower bound for $q$
which is only a power of $r$. However we have been unable to eliminate
entirely the need for a lower bound on $q$, even for $n$ as large as
$2r^2$.

The case $r=8$ is of relevance to the problem of $p$-adic zeros of
quartic forms. The author \cite{hb} has shown that if $p\not=2,5$, any quartic
form over $\Q_p$ in $n$ variables has a non-trivial $p$-adic zero,
providing that any system of 16 linear forms and 8 quadratic forms
also has a non-trivial zero. Our results therefore have the following
corollary.
\begin{corollary}\label{c3}
A quartic form over $\Q_p$ in at least $49$ variables has a
non-trivial $p$-adic zero providing that $p\ge 271919$.
\end{corollary}

Our proofs use a $\cp$-adic minimization technique, for which see Birch
and Lewis \cite[Lemma 12]{BL}.  Let $F$ be the residue  field.  
Then, as in \cite[\S\S 3 \& 4]{BL}, it suffices to prove our
theorem for ``minimized'' systems of forms $q^{(i)}$.  Such forms will 
have $\cp$-adic integer coefficients, and we write
$Q^{(i)}(x_1,\ldots,x_n)\in F[x_1,\ldots,x_n]$ for their reductions
in $F$.  In view
of Hensel's Lemma it will suffice to find a non-singular zero in $F^n$
for the system $Q^{(i)}=0$.  The minimization process ensures that the
forms $Q^{(i)}$ will satisfy a key
condition, given by (2) of \cite[Lemma 12]{BL}.  We proceed to 
explain this condition.  

Suppose $S^{(1)},\ldots,S^{(s)}$ are linearly independent forms taken from the
$F$-pencil generated by the $Q^{(i)}$.  Suppose further that, after a
linear change of variables, the forms
\[ S^{(i)}(0,\ldots,0,x_{w+1},\ldots,x_m)\;\;\; (1\le i\le s)\]
all vanish identically.  Then if the original system $q^{(i)}$ was
minimized, part (2) of \cite[Lemma 12]{BL} tells us that
\begin{equation}\label{cond}
w\ge \frac{sn}{2r}.
\end{equation}
In particular, if $n>4r$ we must have $w>2s$.
As an example of the minimization condition (\ref{cond}), take $n>4r$ 
and $s=1$, whence we deduce that $w\ge 3$.  Thus no form $S$ in the
pencil can be annihilated by setting 2 variables equal to zero.
In particular, if there were any form $S$ in the $F$-pencil which had rank at
most 2 we could express it as a function of $x_1$ and $x_2$ only,
allowing $w=2$, and thereby giving a contradiction.  Indeed if there 
were a form of rank 3 it could be
written as $S(x_1,x_2,x_3)$, and by Chevalley's Theorem we could take
$S(0,0,1)=0$, which again permits $w=2$.  We therefore conclude
that if $n>4r$ the condition (\ref{cond}) implies that very non-zero form in
the $F$-pencil has rank at least 4.

We can now focus on systems $Q^{(i)}$ over the finite field $F$.  As
noted above, it suffices to find a non-singular zero, given the key
minimization condition (\ref{cond}).  This will be done by a counting
argument, in which we first give a lower bound estimate for the total 
number of solutions to the system $Q^{(i)}=0$, and then give an upper
bound on the number of singular solutions.  Here a major 
r\^{o}le will be played by singular forms in
the $F$-pencil generated by the $Q^{(i)}$. We will therefore be forced
to consider how many forms of a given rank the pencil can contain, and
this problem is the key point in the proof. Our treatment will use
some algebraic geometry ultimately motivated by the work of 
Davenport \cite[\S 2]{D16},
and it is at this point that the minimization condition (\ref{cond}) is
applied.
\bigskip

{\bf Acknowledgments.} Part of this work was carried out at the 
Hausdorff Institute for Mathematics in Bonn, during the
Trimestre on Diophantine Equations .  The hospitality and financial 
support of the institute is gratefully acknowledged.

The material on quadratic forms in characteristic 2 owes a great deal
to numerous conversations with Damiano Testa, to whom the author is
delighted to express his gratitude.

\section{Geometric Considerations}

In discussing the geometry of our system of quadratic forms we
shall work over the algebraic closure $\Fb$.  Thus
when we speak of a point on a variety $V$, we shall mean an
$\Fb$-point, unless we explicitly write $V(F)$. We shall take special
care to include the case in which $F$ is dyadic.  We write $\chi(F)$
for the characteristic of $F$.
To begin with we will not assume that condition (\ref{cond}) holds.  

We start by attaching a symmetric $n\times n$ matrix 
$M^{(i)}$ to each form $Q^{(i)}$.  In general, if
\[Q(x_1,\ldots,x_n)=\sum_{1\le i\le j\le n}a_{ij}x_ix_j,\]
then the associated matrix will have entries
\begin{equation}\label{MM}
M_{ij}=\left\{\begin{array}{cc} a_{ij},& i<j,\\ 2a_{ii}, & i=j,\\
a_{ji}, & i>j. \end{array}\right.
\end{equation}
When $\chi(F)\not=2$ this corresponds to the usual
definition.  For $\chi(F)=2$ the matrix $M$ is skew-symmetric, and
always has even rank.

By the rank of a quadratic form $Q$ we mean the minimal $r$ such that
there is a form $Q'$ over $F$, in $r$ variables, and linear forms
\[L_1(x_1,\ldots,x_n), \ldots,L_r(x_1,\ldots,x_n)\]
over $F$ for which 
$Q(x_1,\ldots,x_n)=Q'(L_1,\ldots,L_r)$. It is not hard to show that
the rank of a form is independent of the field over which one works.
When $\chi(F)\not=2$ one has
$\rank(Q)=\rank(M)$, but this is not true in general if $\chi(F)=2$.  
However we always have
\[\rank(M)=2[\rank(Q)/2]\]
for dyadic fields.

When $\chi(F)\not=2$ the condition $\rank(Q)\le R$ is equivalent to the
vanishing of all the $(R+1)\times(R+1)$ minors of $M$.  When $\chi(F)=2$
and $R$ is odd we have $\rank(Q)\le R$ if and only if
$\rank(M)\le R-1$.  Hence in this case a necessary and sufficient
condition is that the $R\times R$ minors of $M$ all vanish.  When 
$\chi(F)\not=2$
and $R$ is even the picture is slightly more
complicated. A necessary and sufficient condition for the rank of $Q$
to be at most $R$ is that $\rank(M)\le R$ and that if $\rank(M)=R$
then $Q$ should vanish on a set of generators for 
the null space of $M$.  However, if
$\rank(M)=R$ then the null space is generated by vectors
$\b{v}_1,\ldots,\b{v}_{n-R}$, whose
components are $R\times R$ minors of $M$, while if $\rank(M)<R$ these
vectors will vanish.  It follows that if $\chi(F)=2$ and
$R$ is even then $\rank(Q)\le R$ if and only if $\rank(M)\le R$ and
$Q(\b{v}_i)=0$ for $i\le n-R$.  Thus in each case there is a set of
polynomial conditions on the coefficients of $Q$, which determines
whether or not $\rank(Q)\le R$.  If we now define
\begin{equation}\label{VR}
V_R=\{[\bu]\in\Proj^{r-1}: \rank\{\sum_{i=1}^r u_iQ^{(i)}\}\le R\}
\end{equation}
it follows that $V_R$ is an algebraic set.
We have shown that that these polynomial conditions defining $V_R$ 
are of degree at most $R+1$ in $\b{u}$ unless $\chi(F)=2$ and $R$
is even, in which case they have degree $2R+1$.  In the final section
of this paper we will establish the following improvement.
\begin{lemma}\label{app}
When $F$ is a finite field with
$\chi(F)=2$ and $R$ is even there is a set of forms of degree
$R+1$ in the coefficients of the quadratic form $Q$ which vanish if
and only if $\rank(Q)\le R$.
\end{lemma}

Suppose that we have a point
$[\bu_0]$ which lies in $V_R(F)$ but not in $V_{R-1}$, where 
we conventionally write
$V_{-1}=\emptyset$. Then $[\bu_0]$ will belong to some component $W$,
say, of $V_R$.  We proceed to bound the dimension of $W$.

Let $k=n-R$ and let $G$ be the Grassmanian of $(k-1)$-dimensional 
linear spaces $L\subseteq\Proj^{n-1}$. Then $\rank(Q)\le R$ if and 
only if there is an $L\in G$ such that $M\x=\b{0}$ and
$Q(\x)=0$ for all $[\x]\in L$.  We use the notation $ML=0$ and
$Q(L)=0$ for these latter conditions. If $\rank(Q)=R$ the space $L$
will be unique, and will be defined over $F$.  If $N$ is the vector
space corresponding to $L$, so that $\dim(N)=k$, we say that $N$ is
the null space for $Q$.  

Let
\[J=\{([\bu],L)\in W\times G:\,(\sum_1^r u_iM^{(i)})L=0,\,
(\sum_1^r u_iQ^{(i)})(L)=0\}.\]
When we project
from $J$ to $W$ the fibre above any point is non-empty, whence 
$\dim(J)\ge\dim(W)$.

It is now convenient to change the basis for the $F$-pencil generated
by the forms $Q^{(i)}$ so that $\bu_0$ becomes $(1,0,\ldots,0)$.  We
then put $Q=Q^{(1)}$, so that $Q$ has rank exactly $R$.  Let $N$ be
the null space for $Q$, and 
make a linear change of variables so that $N$ is
generated by the first $k$ unit vectors $\b{e}_1,\ldots,\b{e}_k$.  We
would like to examine the tangent space of $J$ at $([\bu_0],L_0)$, where
$L_0$ is the projective linear space corresponding to $N$.  This tangent
space is most readily identified by switching to the affine setting.
We therefore define 
\[V=\{\b{v}=(v_2,\ldots,v_r)\in\A^{r-1}:[(1,\b{v})]\in W\}\]
and
\[Y=\{\b{y}\in\A^n: y_1=\ldots=y_k=0\}.\]
Notice that $\b{0}\in V$ and that $\dim(V)=\dim(W)$.

We now consider the algebraic set $Z\subseteq V\times Y^k$ specified
by the condition that $v\in V$, along with the equations
\[\{M+\sum_{i=2}^rv_iM^{(i)}\}(\b{e}_j+\b{y}_j)=\b{0},\;\;\; (1\le j\le k)\]
and
\[\{Q+\sum_{i=2}^rv_iQ^{(i)}\}(\b{e}_j+\b{y}_j)=0,\;\;\; (1\le j\le k).\]
Thus we have $nk+k$ equations, in addition to the condition $v\in
V$. Note that our equations imply that $\{Q+\sum_
iv_iQ^{(i)}\}(\b{w})=0$ for any $\b{w}$ in the span of the
vectors $\b{e}_j+\b{y}_j$.  Thus
$Z$ is an affine version of $J$, with the linear space $L$
corresponding to the vector space generated by $\b{e}_j+\b{y}_j$ for
$1\le j\le k$. In particular it follows that $\dim(Z)=\dim(J)\ge\dim(W)$.

One can now calculate the tangent
space $T=\mathbb{T}(Z,(\b{0},\ldots,\b{0}))$.  One finds that $T$ is
the set of $(\b{v},\b{y}_1,\ldots,\b{y}_{k})\in
\mathbb{T}(V,\b{0})\times Y^k$ which satisfy the equations
\begin{equation}\label{rel1}
\{\sum_{i=2}^rv_iM^{(i)}\}\b{e}_j+M\b{y}_j=\b{0},\;\;\; (1\le j\le
k)
\end{equation}
and
\[\{\sum_{i=2}^rv_iQ^{(i)}\}(\b{e}_j)+\b{y}_j^T\nabla Q(\b{e}_j)=0,
\;\;\; (1\le j\le k).\]
However we have $\nabla Q(\b{e}_j)=M\b{e}_j=\b{0}$, so that the second
set of conditions reduce to
\begin{equation}\label{rel2}
\{\sum_{i=2}^rv_iQ^{(i)}\}(\b{e}_j)=0,\;\;\; (1\le j\le k).
\end{equation}
If $(\b{v},\b{y}_1,\ldots,\b{y}_k)\in T$ we may pre-multiply the relation
(\ref{rel1}) by $\b{e}_h^T$ for any $h\le k$ and use the fact that
$\b{e}_h^TM=\b{0}^T$ to deduce that
\begin{equation}\label{st}
\b{e}_h^T\{\sum_{i=2}^rv_iM^{(i)}\}\b{e}_j=0,\;\;\;(1\le j,h\le k).
\end{equation}
The two conditions (\ref{rel2}) and (\ref{st}) now imply that
\begin{equation}\label{g}
\{\sum_{i=2}^rv_iQ^{(i)}\}(\x)=0\;\;\mbox{for all}\;\; \x\in N.
\end{equation}

Let $\pi:T\rightarrow \mathbb{T}(V,\b{0})$ be the natural 
projection.  Then the 
relation (\ref{g}) holds for any $\b{v}\in\pi(T)$.  However $\pi$ is a
linear map between vector spaces, and 
\[{\rm Ker}(\pi)=\{(\b{0},\b{y}_1,\ldots,\b{y}_{k})\in\{\b{0}\}\times
Y^k:\,M\b{y}_j=\b{0},\;\;\; (1\le j\le k)\}.\]
When $\chi(F)\not=2$ the matrix $M$ has null space $N$, so
that we must have $\b{y}_j=\b{0}$ for all $j$, whence ${\rm Ker}(\pi)$
is trivial.  When $\chi(F)=2$ the matrix $M$ will have null space $N$ only when
$R$ is even.  Thus, in the dyadic case we now require $R$ to be even.
Under this assumption we will have $\dim(\pi(T))=\dim(T)$, whence
\[\dim(\pi(T))=\dim(T)\ge\dim(Z)=\dim(J)\ge\dim(W),\]
since the tangent space of $Z$ at any point has dimension at least as
large as $Z$ itself.

Since $Q^{(1)}(\x)=Q(\x)=0$ for all $\x\in N$ we now deduce that there is
a linear space of quadratic forms in the $\overline{F}$-pencil, with
dimension at least $1+\dim(W)$, all vanishing on the space $N$.  However
$N$ is defined over $F$ itself, whence
\[\{\bu\in \A^r:\{\sum_1^ru_iQ^{(i)}\}(\x)=0\;\,
\mbox{for all}\;\x\in N\}\]
is also defined over $F$.  We therefore draw the following conclusion.

\begin{lemma}\label{fund0}
Let $V_R$ be the variety {\rm (\ref{VR})}.  Suppose either that 
$\chi(F)\not=2$, or that $\chi(F)=2$ and that $R$ is even.  
Suppose further that we have a point 
$\bu\in F^r$ for which the form
\begin{equation}\label{M}
Q=\sum_{i=1}^r u_{i}Q^{(i)},
\end{equation} 
has rank $R$ and null-space $N$, and 
such that $[\bu]$ belongs to 
an irreducible component $W$ of $V_R$.  Then there are at least
$1+\dim(W)$ linearly independent quadratic forms $S^{(i)}$ in the
$F$-pencil {\rm (\ref{M})}, all of which vanish on the $F$-vector space
$N$ of codimension $R$ in $F^n$.
\end{lemma}

To handle the case in which $\chi(F)=2$ and $R$ is odd we need to make
a small modification of the previous argument.  We keep the same
notation as before, but in addition to the
null space $N$ of $Q$ we must now consider the null space $N_0$ of $M$.
In the previous situation these coincided but now $N$ is strictly
contained in $N_0$. If we now write $G_0$ for the Grassmannian of 
$k$-dimensional linear subspaces of $\Proj^{n-1}$ then $N$ and $N_0$ will
corresponds to some pair of linear spaces $L\in G$ and $L_0\in G_0$, 
with $L\subset L_0$.  We now define
\begin{eqnarray*}
\lefteqn{J_0=\{([\bu],L,L_0)\in W\times G\times G_0:\, L\subset L_0,}
\hspace{3cm}\\
&&(\sum_1^r u_iM^{(i)})L_0=0,\,(\sum_1^r u_iQ^{(i)})(L)=0\}.
\end{eqnarray*}
As before, when we project from $J_0$ to $W$, the fibre above any 
point is non-empty, whence  $\dim(J_0)\ge\dim(W)$.

Following the previous analysis we switch to affine coordinates. We
change variables as before, so that $Q=Q^{(1)}$, and so that $N$ and
$N_0$ are generated by $\b{e}_1,\ldots,\b{e}_k$ and 
$\b{e}_1,\ldots,\b{e}_{k+1}$ respectively.  We use the same set $V$ as
before, but take
\[Y=\{\b{y}\in\A^n: y_1=\ldots=y_{k+1}=0\}.\]
This time we define a set $Z_0\subseteq V\times Y^{k+1}$ specified
by the condition that $v\in V$, along with the equations
\[\{M+\sum_{i=2}^rv_iM^{(i)}\}(\b{e}_j+\b{y}_j)=\b{0},\;\;\; (1\le j\le k+1)\]
and
\[\{\sum_{i=2}^rv_iQ^{(i)}\}(\b{e}_j+\b{y}_j)=0,\;\;\; (1\le j\le k).\]
Again we note that $Z_0$ is an affine version of $J_0$, whence 
$\dim(Z_0)=\dim(J_0)\ge\dim(W)$.

The tangent space $T_0=\mathbb{T}(Z_0,(\b{0},\ldots,\b{0}))$ is
the set of $(\b{v},\b{y}_1,\ldots,\b{y}_{k+1})\in
\mathbb{T}(V,\b{0})\times Y^{k+1}$ which satisfy the equations
\[\{\sum_{i=2}^rv_iM^{(i)}\}\b{e}_j+M\b{y}_j=\b{0},\;\;\; (1\le j\le k+1)\]
and
\[\{\sum_{i=2}^rv_iQ^{(i)}\}(\b{e}_j)=0,\;\;\; (1\le j\le k).\]
As before, these imply that
\[\{\sum_{i=2}^rv_iQ^{(i)}\}(\x)=0\;\;\mbox{for all}\;\; \x\in N.\]
If $\pi_0:T_0\rightarrow \mathbb{T}(V,\b{0})$ is the natural 
projection then the above relation holds for any $\b{v}\in\pi(T_0)$.  
However 
\[{\rm Ker}(\pi_0)=\{(\b{0},\b{y}_1,\ldots,\b{y}_{k+1})\in\{\b{0}\}\times
Y^{k+1}:\,M\b{y}_j=\b{0},\;\;\; (1\le j\le k+1)\}.\]
Since $M$ has null space $N_0$ we must have $\b{y}_j=\b{0}$ for all
$j$, whence ${\rm Ker}(\pi_0)$ is trivial. We may now complete the
argument as before, leading to the following conclusion.

\begin{lemma}\label{fund1}
Let $V_R$ be the variety {\rm (\ref{VR})}.  Suppose that $\chi(F)=2$ 
and that $R$ is odd. Suppose further that we have a point 
$\bu\in F^r$ for which the form
\begin{equation}\label{M0}
Q=\sum_{i=1}^r u_{i}Q^{(i)},
\end{equation} 
has rank $R$ and null-space $N$, and 
such that $[\bu]$ belongs to 
an irreducible component $W$ of $V_R$.  Then there are at least
$1+\dim(W)$ linearly independent quadratic forms $S^{(i)}$ in the
$F$-pencil {\rm (\ref{M0})}, all of which vanish on the $F$-vector space
$N$ of codimension $R$ in $F^n$.
\end{lemma}

If we now assume the fundamental minimization condition (\ref{cond}) then
we may take $n-w=\dim(N)$, so that 
\[R=n-\dim(N)=w\geq \frac{n}{2r}(1+\dim(W)),\]
and therefore $1+\dim(W)\le 2rR/n$.  

\begin{lemma}\label{fund}
Suppose that {\rm (\ref{cond})} holds. Let $V_R$ be the variety 
{\rm (\ref{VR})}.  Then any point $[\bu]\in \mathbb{P}^{r-1}(F)$ for
which the form {\rm (\ref{M})} has rank $R$ will belong to an irreducible
component $W$ of $V_R$ having $1+\dim(W)\le 2rR/n$.
\end{lemma}

This lemma is the most novel part of our argument.  Notice that it
tells us nothing about those components $W$ of $V_R$ which do not
contain a point defined over $F$, or for which the only such points
are in the subvariety $V_{R-1}$.
\bigskip

We next estimate how many points can lie in each component $W$.
\begin{lemma}\label{count}
Suppose that $V\subseteq \mathbb{A}^r$ is an 
algebraic set of dimension $w$ and degree $d$.  Then
\[\# V(F)\le dq^w,\]
where $q=\# F$.
\end{lemma}
This is a relatively standard result, proved along the lines given by
Browning and the author \cite[page 91]{bhb}. We use induction on $w$,
the case $w=0$ being trivial.  Clearly we can assume that $V$ is
absolutely irreducible, by additivity of the degree.  When $w\ge 1$
there is always at least one index $i$ such that $V$ intersects the
hyperplane $u_i=\alpha$ properly for every $\alpha\in\Fb$. (If this 
were not the case, then $V$ must be contained in a hyperplane
$u_i=\alpha_i$ for each index $i$, so that $V$ could contain at most
the single point $(\alpha_1,\ldots,\alpha_r)$.)  Fixing a suitable
index $i$ we conclude that
\[\# V(F)\le \sum_{\alpha\in F} \#(V\cap\{u_i=\alpha\}).\]
Since $V\cap\{u_i=\alpha\}$ has dimension at most $w-1$ and degree at
most $d$ we may use the induction hypothesis to conclude that
\[\#(V\cap\{u_i=\alpha\})\le dq^{w-1},\]
whence the required induction bound follows.
\bigskip

In order to estimate the contribution from all the relevant components
$W$ of $V_R$ we will need information on their degrees as well as
their dimensions, and for this we use the following result.
\begin{lemma}\label{deg}
Let $V\subseteq \mathbb{A}^r$ be an algebraic set defined by the
vanishing of polynomials $f_1,\ldots,f_N$ each having total degree
at most $d$.  Suppose that $V$ decomposes into irreducible components as
$V=\cup_{i=1}^I V_i$.  Then
\[\sum_{i=1}^I\deg(V_i)d^{\dim(V_i)}\le d^r.\]
\end{lemma}
This is proved by induction on $N$, the case $N=1$ being trivial.  We
proceed to assume that the result holds for the case $N$, and prove it
for the case $N+1$.  Let us write $H=\{f_{N+1}=0\}$ for convenience,
and suppose that $V_i\cap H$ decomposes into irreducible components as
$\cup_{j=1}^{J(i)} V_{ij}$.  We claim that
\begin{equation}\label{cl}
\sum_{j=1}^{J(i)}\deg(V_{ij})d^{\dim(V_{ij})}\le
\deg(V_i)d^{\dim(V_i)}.
\end{equation}
Once this is established we will have
\[\sum_{i=1}^I\sum_{j=1}^{J(i)}\deg(V_{ij})d^{\dim(V_{ij})}\le
\sum_{i=1}^I\deg(V_i)d^{\dim(V_i)}\le d^r,\]
by the induction hypothesis.  We will therefore have completed the
induction step.

To prove the statement (\ref{cl}) we factor $f_{N+1}$ into absolutely
irreducible polynomials $f_{N+1}=g_1\ldots g_M$, say, and write
$H_k=\{g_k=0\}$. If there is any index $k$ such that $V_i\subseteq
H_k$ then $V_i\subseteq H$, whence $V_i\cap H=V_i$ is already
irreducible and (\ref{cl}) is trivial.  On the other hand, if $V_i$
and $H_k$ intersect properly for every $k$ then $V_i\cap H_k$ is a
union of components $V_{ij}$ for $j$ in some set
$S(k)\subseteq\{1,\ldots,J(i)\}$, with $\dim(V_{ij})=\dim(V_i)-1$ and
\[\sum_{j\in S(k)}\deg(V_{ij})\leq \deg(V_i)\deg(g_k)\]
by B\'{e}zout's Theorem.  Summing over $k$ then yields
\[\sum_{j=1}^{J(i)}\deg(V_{ij})\le\deg(V_i)d,\]
and (\ref{cl}) follows in this case too.  This completes the proof of
Lemma \ref{deg}.

We now combine Lemmas \ref{fund}, \ref{count} and \ref{deg} to produce
the following result.
\begin{lemma}\label{sum}
Suppose that the quadratic forms $q^{(i)}$ form a minimized system.
Then the number $N(R)$ of quadratic forms {\rm (\ref{M})} of rank $R$, 
with $\bu\in F^r$, satisfies
\[N(R)\le (\frac{q}{R+1})^{[2rR/n]}( R+1)^r\]
whenever $q\ge  R+1$.  Moreover  any non-zero form in the $F$-pencil
has rank at least $2(\lceil n/2r\rceil-1)$.
\end{lemma}

Suppose that $V_R$ is a union 
\[V_R=\bigcup_1^I W_i\]
of irreducible
components, and that the points $[\bu]\in V_R(F)$
lie in components $W_1,\ldots, W_L$.  Then, applying Lemma \ref{count}
to the affine cone over each $W_i$, we find that
\[N(R)\le \sum_{i=1}^L\deg(W_i)q^{1+\dim(W_i)}.\]
However $V_R$ is defined by equations of degree at most $
R+1=d$, say, whence Lemma \ref{deg} yields
\[\sum_{i=1}^L\deg(W_i)( R+1)^{1+\dim(W_i)}\le 
\sum_{i=1}^I\deg(W_i)( R+1)^{1+\dim(W_i)}\le ( R+1)^r.\]
However Lemma \ref{fund} shows that $1+\dim(W_i)\le [2rR/n]$ for $i\le
L$, so that if $q\ge  R+1$ we will have
\begin{eqnarray*}
N(R)&\le&\sum_{i=1}^L\deg(W_i)( R+1)^{1+\dim(W_i)}
(\frac{q}{ R+1})^{1+\dim(W_i)}\\
&\le& (\frac{q}{ R+1})^{[2rR/n]}
\sum_{i=1}^L\deg(W_i)( R+1)^{1+\dim(W_i)}\\
&\le& (\frac{q}{ R+1})^{[2rR/n]}( R+1)^r
\end{eqnarray*}
as required. 

For the final observation we extend the remark made in \S 1, in 
connection with the condition (\ref{cond}).  Any form of rank $R$ over
$F$ will vanish on a vector space of codimension $(R+1)/2$, if $R$ is
odd, or of codimension $(R+2)/2$ if $R$ is even.  We may therefore
take $w=1+[R/2]$ and deduce that $1+[R/2]\ge n/2r$, which gives the
required lower bound on $R$. Note that this argument uses only the
minimization condition, and does not require either Lemmas \ref{fund0},
\ref{fund1} or \ref{fund}.

\section{Counting Zeros}

We begin by considering zeros of a system of quadratic forms 
\[S^{(i)}(x_1,\ldots,x_k)\in F[x_1,\ldots,x_k],\;\;\;(1\le i\le I).\]
Consider the set
\[A=\{(\bu,\x)\in F^I\times F^k: \sum_{i=1}^I u_i
S^{(i)}(x_1,\ldots,x_k)=0\}.\]
We shall count elements of $A$ in two ways.  Firstly we consider how
many choices of $\bu$ correspond to each $\x$.  If $S^{(i)}(\x)=0$ for
each index $i$ then there are $q^I$ possible vectors $\bu$, and
otherwise $q^{I-1}$ choices.  Hence if the system $S^{(i)}(\x)=0$ has
$N$ zeros in total we will have
\[\# A=q^IN+q^{I-1}(q^k-N).\]
Alternatively we can count elements of $A$ according to the value
$\bu$.  Here we write
\[N(\bu)=\#\{\x\in F^k: \sum_{i=1}^I u_i
S^{(i)}(x_1,\ldots,x_k)=0\},\]
whence
\[\# A=\sum_{\lb{u}}N(\b{u}).\]
We therefore deduce that
\begin{eqnarray*}
N&=&
\frac{1}{q^{I-1}(q-1)}\left\{-q^{I+k-1}+\sum_{\lb{u}}N(\b{u})\right\}\\
&=&\frac{1}{q^{I-1}(q-1)}\left\{\sum_{\lb{u}}(N(\b{u})-q^{k-1})\right\}\\
&=&q^{k-I}+
\frac{1}{q^{I-1}(q-1)}\left\{\sum_{\lb{u}\not=0}(N(\b{u})-q^{k-1})\right\},
\end{eqnarray*}
since $N(\b{0})=q^k$.

We proceed to consider the number $N(S)$ of zeros of a single quadratic
form $S(x_1,\ldots,x_k)$.  If $\rank(S)=0S$, then there
are trivially $q^k$ zeros, and if $S$ has rank one there are $q^{k-1}$
zeros.  For rank 2 there will be $(2q-1)q^{k-2}$ zeros if $S$ factors over $F$
and $q^{k-2}$ zeros otherwise.  For larger ranks there will be at least one
non-singular zero, by
Chevalley's Theorem, and a linear change of variable will allow us to write
$S$ in the shape 
\[S(x_1,\ldots,x_k)=x_1x_2+S'(x_3,\ldots,x_k).\]
One then finds that there are $2q-1$ possibilities for $(x_1,x_2)$ if
$S'=0$ and $(q-1)$ choices otherwise, so that
$N(S)=qN(S')+(q-1)q^{k-2}$.  An easy induction on $k$ now shows that
$N(S)=q^{k-1}$ whenever $S$ has odd rank, and that
\[|N(S)-q^{k-1}|=(1-q^{-1})q^{k-R/2}\]
whenever $S$ has even rank $R$.

We may therefore conclude as follows.
\begin{lemma}\label{count2}
Suppose we have a system of quadratic forms 
\[S^{(i)}(x_1,\ldots,x_k)\in F[x_1,\ldots,x_k],\;\;\;(1\le i\le I)\]
with $N$ zeros over $F$.  Write $N_R$ for the number of vectors
$\bu\in F^I$ for which
\begin{equation}\label{comb}
\sum_{i=1}^I u_i S^{(i)}(x_1,\ldots,x_k)
\end{equation}
has rank $R$, and assume that such a linear combination vanishes only
for $\bu=\b{0}$.  Then
\[|N-q^{k-I}|\le\sum_{1\le t\le k/2}q^{k-I-t}N_{2t}.\]
\end{lemma}

We may now apply Lemma \ref{count2} to count non-singular zeros of the
system
\begin{equation}\label{sys}
Q^{(1)}(x_1,\ldots,x_n),\ldots,Q^{(r)}(x_1,\ldots,x_n)
\end{equation}
arising from a minimized system $q^{(1)},\ldots,q^{(r)}$.  In view of
Lemmas \ref{count} and \ref{sum} the total number $N$ of common zeros
satisfies
\begin{equation}\label{N}
N\ge q^{n-r}\left\{1-\sum_{\lceil n/2r\rceil-1\le t\le n/2}q^{-t}
(\frac{q}{2 t+1})^{[4rt/n]}(2 t+1)^r\right\}
\end{equation}
providing that $q>  n\ge 4r+1$.  This latter 
condition is enough to
ensure that $q\ge 2 t+1$ whenever $t\le n/2$. 
Note that if a non-trivial linear combination (\ref{comb}) were to
vanish we would be able to take $s=1, w=0$ in (\ref{cond}), which is
impossible. We remark that the sum in (\ref{N}) is $O_{r,n}(q^{-1})$ as soon as
$n>4r$, and indeed we will have $N\sim q^{n-r}$ as
$q\rightarrow\infty$, for such $n$.  This is the behaviour we would
have if the variety defined by $q^{(1)}=\ldots=q^{(r)}=0$ were
absolutely irreducible.  However it is not clear whether the
minimization condition ensures such irreducibility.
\bigskip

We have now to consider singular zeros for the system (\ref{sys}).
Any such zero $\x$ is a singular zero of at least one non-zero
form (\ref{comb}) in
the pencil, $S$ say. Unless $\x=\b{0}$ we may deduce that $S$ is
singular.  We proceed to estimate how many zeros the system (\ref{sys})
has, which are singular zeros of a given form $S$ of the shape (\ref{comb}).
By changing the basis for the pencil we may indeed
assume that $S=Q^{(r)}$.  Suppose that $S$ has rank $\rho<n$.  Then the 
singular zeros of $S$ form a
vector space of dimension $n-\rho=k$, say, which we may take to be
\[\{(x_1,\ldots,x_k,0,\ldots,0)\},\]
after a suitable change of variable. It follows then that our problem
is to count zeros of the new system
\[S^{(1)}(x_1,\ldots,x_k),\ldots,S^{(r-1)}(x_1,\ldots,x_k),\]
where
\[S^{(i)}(x_1,\ldots,x_k)=Q^{(i)}(x_1,\ldots,x_k,0,\ldots,0).\]
According to Lemma \ref{count2} there are at most
\begin{equation}\label{B}
q^{k-(r-1)}\{\sum_{0\le t\le k/2}q^{-t}N_{2t}\}
\end{equation}
such zeros, where $N_R$ is the number of linear combinations 
\begin{equation}\label{sys1}
\sum_{i=1}^{r-1} u_i S^{(i)}(x_1,\ldots,x_k)
\end{equation}
which have rank $R$.

To estimate $N_R$ we will use Lemmas \ref{fund0} and \ref{fund1} 
in combination with
Lemmas \ref{count} and \ref{deg}. If $R=2t$ and $W\subseteq\mathbb{P}^{r-2}$ 
is an irreducible component of the variety of vectors counted by
$N_R$, then Lemmas \ref{fund0} and \ref{fund1} show that we have at 
least $1+\dim(W)$ linearly independent forms
from the pencil (\ref{sys1}) which vanish simultaneously on a vector space
$X\subseteq F^k$ of codimension $R$.  By extending these to forms on 
$F^r$ we obtain
$1+\dim(W)$ linearly independent forms from the pencil 
\[\sum_{i=1}^{r-1} u_i Q^{(i)}(x_1,\ldots,x_n)\]
which vanish simultaneously on 
\[\tilde{X}=\{(x_1,\ldots,x_k,0,\ldots,0)\in F^n:\,
(x_1,\ldots,x_k)\in X\}.\]
However $Q^{(r)}$ also vanishes on $\tilde{X}$, whence the minimization
condition (\ref{cond}) yields 
\[n-\dim(\tilde{X})\geq\frac{(2+\dim(W))n}{2r}.\]
Since $\dim(\tilde{X})=\dim(X)=k-R$ we deduce that
\begin{equation}\label{wc}
\dim(W)\le \frac{2r(n-k+R)}{n}-2.
\end{equation}
This allows us to use Lemmas \ref{count} and \ref{deg} to conclude that
\[N_R\le (\frac{q}{ R+1})^{[2r(n-k+R)/n]-1}( R+1)^{r-1},\]
for $q\ge  R+1$, as in the proof of Lemma \ref{sum}.  

Since $k=n-\rho$ we now find from
(\ref{B}) that the number of zeros of (\ref{sys}) which are singular
for a particular $S$ of rank $\rho$ is at most
\begin{eqnarray*}
\lefteqn{q^{n-\rho-r+1}\{\sum_{0\le t\le (n-\rho)/2}q^{-t}
(\frac{q}{2 t+1})^{[2r(\rho+2t)/n]-1}(2t+1)^{r-1}\}}\hspace{1cm}\\
&=&q^{n-\rho-r}\{\sum_{0\le t\le (n-\rho)/2}q^{-t}
(\frac{q}{2 t+1})^{[2r(\rho+2t)/n]}(2t+1)^{r}\}.
\end{eqnarray*}
To estimate the total number of singular zeros of (\ref{sys}) we must
sum this over all singular forms $S$, and allow for the trivial
singular zero $\x=\b{0}$.  Although Lemma \ref{sum}
estimates the number of singular forms of given rank, for our present
purposes scalar multiples of a given form $S$ produce the same
singular zeros.  Hence it suffices to count only one form $S$ from
each set of scalar multiples.  Thus Lemma \ref{sum} shows that the
total number of non-trivial 
singular zeros for the system (\ref{sys}) is at most
\[\frac{q^{n-r}}{q-1}\sum_{\rho=2(\lceil n/2r\rceil-1)}^{n-1}
(\frac{q}{\rho+1})^{[2r\rho/n]}
\frac{(\rho+1)^r}{q^{\rho}}\sum_{0\le t\le (n-\rho)/2}
(\frac{q}{2 t+1})^{[2r(\rho+2t)/n]}\frac{(2 t+1)^{r}}{q^t}\]
for $q> n$.  Note that this latter condition will ensure that $q\ge
2 t+1$ and that $q\ge \rho+1$.  After allowing for $\x=\b{0}$ it 
now follows that the total number of 
non-singular zeros for the system
(\ref{sys}) is at least $q^{n-r}(1-\sigma_1-\sigma_2)$ with $\sigma_1$
and $\sigma_2$ as in the theorem, and the sufficiency of the condition
$\sigma_1+\sigma_2<1$ follows.

\section{Completion of the Proofs}

We begin by examining the special case $n=4r+1$.  With this value of
$n$ we have $[4rt/n]=t-1$ for $2\le t\le n/2$, whence
\[\sigma_1=q^{-3r-1}+q^{-1}\sum_{2\le t\le 2r}(2t+1)^{r-t+1}.\]
Similarly we will have $[2r\rho/n]=(\rho-1)/2$ and
$[2r(\rho+2t)/n]=t+(\rho-1)/2$ if $\rho$ is odd, while
$[2r\rho/n]=\rho/2-1$ and $[2r(\rho+2t)/n]=t+\rho/2-1$ if $\rho$ is
even.  Thus
\begin{eqnarray*}
\sigma_2&=&\frac{1}{q-1}\left\{q^{-1}
\sum_{\nu=2}^{2r-1}\,\sum_{0\le t\le 2r-\nu}
(2\nu+2)^{r-\nu}(2t+1)^{r-t-\nu}\right.\\
&&\hspace{2cm}\mbox{}+
\left. q^{-2}\sum_{\nu=2}^{2r}\,\sum_{0\le t\le 2r-\nu}
(2\nu+1)^{r-\nu+1}(2t+1)^{r-t-\nu+1}\right\}.
\end{eqnarray*}

In the case $r=3$ we calculate that
\[\sigma_1=q^{-10}+(32.11\ldots)q^{-1}\]
and
\[\sigma_2=(14.72\ldots)q^{-1}(q-1)^{-1}+(145.68\ldots)q^{-2}(q-1)^{-1},\]
whence $q\ge 37$ is admissible.  The other values for $r=4$ and
8 are calculated similarly. 

To prove the general bound it now suffices to assume that $r\ge 5$.
We observe that $(2t+1)^{r-t+1}\le (2r)^{r-1}$ for $2\le t\le
r-1$, while $(2t+1)^{r-t+1}\le 4r+1$ for $r\le t\le 2r$.  It follows that
\[\sum_{2\le t\le 2r}(2t+1)^{r-t+1}\le (r-2)(2r)^{r-1}+(r+1)(4r+1)
\le (r-1)(2r)^{r-1}.\]
For $\sigma_2$ we note that $(2\nu+2)^{r-\nu}\le (2r)^{r-2}$ in each
of the cases $2\le\nu\le r-1$ and $r\le\nu\le 2r-1$, and similarly that
$(2t+1)^{r-t-\nu}\le (2r)^{r-2}$ in all cases.  Thus
\[\sum_{\nu=2}^{2r-1}\,\sum_{0\le t\le 2r-\nu}
(2\nu+2)^{r-\nu}(2t+1)^{r-t-\nu}\le (2r)^{2r-2}.\]
In the same way we have $(2\nu+1)^{r-\nu+1}\le (2r+1)^{r-1}$ and
$(2t+1)^{r-t-\nu+1}\le (2r-1)^{r-1}$ in all cases, whence
\[\sum_{\nu=2}^{2r}\,\sum_{0\le t\le 2r-\nu}
(2\nu+1)^{r-\nu+1}(2t+1)^{r-t-\nu+1}\le (2r)^{2r}.\]
The condition $\sigma_1+\sigma_2<1$ is therefore satisfied if
\[q^{-r}+(r-1)(2r)^{r-1}q^{-1}+(2r)^{2r-2}q^{-1}(q-1)^{-1}+
(2r)^{2r}q^{-2}(q-1)^{-1}<1.\]
One now readily verifies that the above inequality holds if
$r\ge 5$ and $q\ge (2r)^r$, as required for the theorem.

We turn now to Corollary \ref{c2}. Since $n\ge r^2+1$ we have $\lceil
n/2r\rceil-1\ge (r-1)/2$. Thus if $\phi=1-4/r$ we have
\[\sigma_1\le q^{-r}+\sum_{t\ge (r-1)/2}q^{-\phi t}(2t+1)^r.\]
In the infinite sum the ratio of the terms for $t+1$ and $t$ is
\[q^{-\phi}(1+\frac{2}{2t+1})^r\le q^{-\phi}(1+\frac{2}{r})^r\le
q^{-\phi}e^2.\]
Moreover, for a real variable $t$ the function $q^{-\phi t}(2t+1)^r$
is decreasing for $t\ge (r-1)/2$, providing only that $q^{\phi}>e^2$.
It follows that the first term in the sum is at most
$q^{-\phi(r-1)/2}r^r$, whence
\begin{equation}\label{E}
\sum_{t\ge (r-1)/2}q^{-\phi t}(2t+1)^r
\le \frac{q^{-\phi(r-1)/2}r^r}{1-q^{-\phi}e^2}
\end{equation}
and
\[\sigma_1\le q^{-r}+\frac{q^{-\phi(r-1)/2}r^r}{1-q^{-\phi}e^2}\]
if $r\ge 5$ and $q^{\phi}> e^{2}$.

Similarly we find that
\[\sigma_2\le\frac{1}{q-1}\left\{\sum_{\rho= r-1}^{\infty}
\sum_{t=0}^{\infty}q^{-\rho\phi-t\phi}(\rho+1)^r(2t+1)^r\right\}.\]
The double sum factors, and the summation over $\rho$ is
\[\sum_{\rho= r-1}^{\infty}q^{-\rho\phi}(\rho+1)^r\le
\frac{q^{-\phi(r-1)}r^r}{1-q^{-\phi}e}\]
by an argument closely analogous to that above. For the $t$-summation
we note that the real variable function $f(\tau)=\tau^r
q^{-\phi\tau/2}$ is maximal at $\tau=2r/(\phi\log q)$, with maximum
value $\{2r/(e\phi\log q)\}^r\le (r/e)^r$ if $q^{\phi}>e^2$. Thus
\[\sum_{0\le t\le (r-2)/2}q^{-\phi t}(2t+1)^r\le
\frac{r}{2}q^{\phi/2}(r/e)^r.\] 
On combining this with (\ref{E}) we deduce that
\[\sigma_2\le\frac{1}{q-1}
\left\{\frac{q^{-\phi(r-1)}r^r}{1-q^{-\phi}e}\right\}
\left\{\frac{r}{2}q^{\phi/2}(r/e)^r+
\frac{q^{-\phi(r-1)/2}r^r}{1-q^{-\phi}e^2}\right\}.\]
Assuming that $q^{\phi}\ge 2e^2$ we conclude that
\begin{eqnarray*}
\sigma_2&\le& q^{-\phi(r-3/2)}r^{2r}\frac{1}{q-1}
\left\{\frac{1}{1-1/2e}\right\}
\left\{\frac{r}{2}e^{-r}+2q^{-\phi r}\right\}\\
&\le& q^{-\phi(r-3/2)}r^{2r}\frac{2}{q}
\{\frac{r}{2}e^{-r}+2e^{-2r}\}\\
&\le& q^{-\phi(r-1/2)}r^{2r}C_r
\end{eqnarray*}
where
\[C_r=\{\frac{r}{2}e^{-r}+2e^{-2r}\}\le 1\] 
for $r\ge 5$.

One may now calculate that $\phi_1+\phi_2<1$ providing that
$q^{\phi}\ge 4r^2(\ge 2e^2)$.  However, the function $(2r)^{1/(r-4)}$
is decreasing for $r\ge 5$, so that
\[(4r^2)^{1/\phi}=(4r^2)\{(2r)^{1/(r-4)}\}^8\le 10^8(4r^2),\]
and Corollary \ref{c2} follows.

\section{Ranks of Quadratic Forms in \newline Characteristic 2}

In this final section we will prove Lemma \ref{app}.  Let $t_{ij}$ be
indeterminates for $1\le i\le j\le n$, and write
$\b{t}=(t_{11},t_{12},\ldots,t_{nn})$. Let
\begin{equation}\label{aqd}
Q_{\lb{t}}(x_1,\ldots,x_n)=\sum_{1\le i\le j\le n}t_{ij}x_ix_j
\end{equation}
be the corresponding quadratic form, considered as a polynomial in
\[\Z[t_{11},t_{12},\ldots,t_{nn},x_1,\ldots,x_n].\]
We associate a  matrix
$U(\b{t})$ to $Q_{\lb{t}}$, with entries
\[U_{ij}=\left\{\begin{array}{cc} t_{ij},& i<j,\\ 2t_{ii}, & i=j,\\
t_{ji}, & i>j. \end{array}\right.\]
If $I,J\subseteq
\{1,\ldots,n\}$ with $\#I=\#J=R+1$ we define $m^*_{I,J}(\b{t})$ to be the
$I,J$ minor of $U$.  This has order $(R+1)\times(R+1)$, and is a form
of degree $R+1$ in the variables $t_{ij}$. If $R$ is even, as we are
supposing, then $m^*_{I,I}(\b{t})$ vanishes modulo 2, since it becomes the
determinant of a skew-symmetric matrix of odd order when we reduce to
$\Z_2$. Thus if we define
\[m_{I,J}(\b{t})=\left\{\begin{array}{cc} m^*_{I,J}(\b{t}), & I\not=J,\\
\frac{1}{2}m^*_{I,I}(\b{t}), & I=J,\end{array}\right.\]
then $m_{I,J}$ will be an integral form in the $t_{ij}$.

We now map the various $m_{IJ}(\b{t})$ to forms $m_{IJ}(\b{t};F)$ 
in $F[t_{11},\ldots,t_{nn}]$, using the obvious homomorphism from 
$\Z[t_{11},\ldots,t_{nn}]$ to $F[t_{11},\ldots,t_{nn}]$. Let
\[Q(x_1,\ldots,x_n)=\sum_{1\le i\le j\le n}a_{ij}x_ix_j\]
be a quadratic form over a finite field $F$ of characteristic 2.
Then we claim
that a necessary and sufficient condition for $Q$ to have rank at most $R$, is
that the forms $m_{I,J}(\b{t};F)$ all vanish at $t_{ij}=q_{ij}$.
This will clearly suffice for Lemma \ref{app}.  It will be convenient
to call this condition on $Q$ the ``Rank Condition''.

Any quadratic form over $F$ can be reduced, via a sequence of
elementary transformations, into a form of the
shape
\[x_1x_2+\ldots+x_{2m-2}x_{2m}+q(x_{2m+1},\ldots,x_{n}),\]
in which 
\[q(x_{2m+1},\ldots,x_{n})=0,\;\; \mbox{or}\;\; x_{2m+1}^2
,\;\; \mbox{or}\;\; x_{2m+1}^2+x_{2m+1}x_{2m+2}+\mu x_{2m+2}^2.   \]
In the third case
$\mu\in F$ is such that $q$ is irreducible over $F$.  The rank of the
form will be $2m$ or $2m+1$ or $2m+2$ respectively.  One can easily
verify by explicit calculation that our claim holds if $Q$ is in one
of these three canonical shapes.

We proceed to show that if forms $Q$ and $Q'$, with coefficients
$q_{ij}$ and $q'_{ij}$ respectively, are related by an elementary
transformation, then $Q$ satisfies the Rank Condition if and only if
$Q'$ does.  This will be sufficient
to complete the proof.  Indeed, since elementary transformations are
invertible, it will be enough to assume that $Q$ satisfies the Rank
Condition and to deduce that $Q'$ does.

Elementary transformations are of three types.  The first kind
interchanges two of the variables $x_i$ and $x_j$, and in this case our
result is trivial, since the forms $m_{I,J}(\b{t};F)$ will merely be
permuted.  The second type of transformation is $S(\lambda)$, say, which 
multiplies $x_1$ by a non-zero scalar $\lambda$.  If we apply $S(v)$,
with an indeterminate $v$, to the quadratic form (\ref{aqd}), then
the forms $m_{I,J}^*(\b{t})$ will be multiplied by appropriate powers of
$v$. It follows that $S(\lambda)$ will
multiply each $m_{I,J}(\b{q};F)$ by a power of $\lambda$.  Hence
again we see that if $Q$ satisfies the Rank Condition then so does $Q'$.

The third type of elementary transformation, which we denote by $T(\lambda)$, 
replaces $x_1$ by $x_1+\lambda x_2$.  The argument here is similar to
that used for $S(\lambda)$.  When $T(v)$ is applied to $Q_{\lb{t}}$
the forms $m_{I,J}^*(\b{t})$ get replaced by linear
combinations of various $m_{K,L}^*(\b{t})$, with coefficients $1,v$ or
$v^2$.  Hence when $T(\lambda)$ is applied to $Q$
the forms $m_{I,J}(\b{q};F)$ get replaced by linear
combinations of various $m_{K,L}(\b{q};F)$, with coefficients 
$1,\lambda$ or $\lambda^2$.  Again it is clear that
if $Q$ satisfies the Rank Condition then so does $Q'$. This completes the proof
of the lemma.

\bigskip
\bigskip

Mathematical Institute,

24--29, St. Giles',

Oxford

OX1 3LB

UK
\bigskip

{\tt rhb@maths.ox.ac.uk}

\end{document}